\documentclass[12pt,twoside,centertags,reqno]{article}
\usepackage{amsfonts,amssymb,amsmath,amsthm,latexsym,amsxtra,amscd}
\usepackage[english]{babel}
\usepackage{verbatim,alltt}
\usepackage{epsfig,epic,eepic}
\usepackage{fancyheadings}
\usepackage{enumerate}
\usepackage{euscript}
\newtheorem{lemma}{Lemma}[section]
\newtheorem{proposition}[lemma]{Proposition}
\newtheorem{theorem}[lemma]{Theorem}
\newtheorem{corollary}[lemma]{Corollary}
\theoremstyle{definition}
\newtheorem{definition}[lemma]{Definition}
\newtheorem{example}[lemma]{Example}

\newtheorem{remark}[lemma]{Remark}
\newtheorem{problem}[lemma]{Problem}
\date{}
\newcommand{\addresses}[3]{}
\newcommand{\emails}[3]{}
\newcommand{\classification}[1]{%
\renewcommand{\thefootnote}{}%
\footnotetext{\mbox{\hspace*{-12pt}
2001 {\it Mathematics Subject Classification.}}
{#1}%
}}
\newcommand{\keywords}[1]{\renewcommand{\thefootnote}{}%
\footnotetext{\mbox{\hspace*{-12pt}
{\it Key words: }}
{#1}%
}}
\title{{\bf Superisolated Surface Singularities}}
\author{ Enrique Artal Bartolo \and Ignacio Luengo \and Alejandro Melle Hern\'andez}
\addresses{Departamento de Matem\'aticas\\
Universidad de Zaragoza\\
Campus Pza. San Francisco s/n\\
E-50009 Zaragoza SPAIN
}
{Departamento de \'Algebra\\
Universidad Complutense\\
Plaza de Ciencias 3\\
E-28040 Madrid SPAIN
}
{Departamento de \'Algebra\\
Universidad Complutense\\
Plaza de Ciencias 3\\
E-28040 Madrid SPAIN
}
\emails{artal@unizar.es
}
{iluengo@mat.ucm.es
}
{amelle@mat.ucm.es
}
\usepackage{hyperref}
\newcommand\bc{{\mathbb C}} \newcommand\ba{{\mathbb A}}
\newcommand\bp{{\mathbb P}} \newcommand\bn{{\mathbb N}}
 \newcommand\bz{{\mathbb Z}}
\newcommand\bq{{\mathbb Q}} 
\newcommand\br{{\mathbb R}}

\newcommand\cld{{\mathcal D}}

  \newcommand\bd{{\mathbb D}}

\DeclareMathOperator\sing{Sing}
\DeclareMathOperator\spec{Sp}

\DeclareMathOperator\mult{mult}

\newcommand\cp{{\mathcal P}}

\DeclareMathOperator\irr{Irr}
\DeclareMathOperator\ltop{top}
\DeclareMathOperator\Top{Top}

\theoremstyle{plain}
\newtheorem{cjt}[lemma]{Conjecture}

\begin{document}
\setcounter{page}{13}
\maketitle
\thispagestyle{empty}
\centerline{\emph{Dedicated to Gert-Martin Greuel on the occasion of his 60th
    birthday.}} 
\classification{%
%    INSERT here MSC Classification
%
32Sxx, 4B05,14H20,14J17}
\renewcommand{\thefootnote}{\arabic{footnote}}%

\begin{abstract}
\noindent
%
%    INSERT here ABSTRACT of your article
%
In this survey, we review part of the  theory of superisolated surface singularities
(SIS)  
and its applications including some new and recent developments. The class of
SIS singularities is, in some sense, the simplest class of germs 
of normal surface singularities. Namely, their tangent cones are reduced
curves and  the geometry and topology of the SIS singularities can be deduced
from them. Thus this class \emph{contains}, in a canonical way, all
the complex projective plane curve theory, which gives a series of nice
examples and counterexamples. They were introduced by I.~Luengo to show the
non-smoothness of the $\mu$-constant stratum and have been used to answer negatively
some other interesting open questions. We review them and the new results on normal
surface singularities whose link are rational homology spheres. 
We also discuss some positive results which have been proved for SIS
singularities.

\end{abstract}

%
%    INSERT here BODY of your article
%

%
%%%%%%%%%%%%%%%%%%%%%%%%%%%%%%%%%%%%%%%%%%%%%%%%%%%%%%%%%%%%%%%%%%%%%%%%%%%%
%      At the end of the article: acknowledgements and literature
%
%\section*{Acknowledgements}
%

%
%     INSERT here acknowledgements
%
\addcontentsline{toc}{section}{Introduction}
\section*{Introduction}

A superisolated surface, SIS for short, singularity $(V,0)\subset (\bc^3,0)$ is a generic
perturbation of the cone over a (singular) reduced projective  plane curve $C$ of degree $d$,
$C=\{f_d(x,y,z)=0\}\subset \bp^2$, by monomials of higher degree. The geometry, resolution 
and topology of $(V,0)$ is determined by the singularities of $C$ and the pair
$(\bp^2,C)$. This provides a canonical way to \emph{embed} the classical and rich theory of complex 
projective plane curves 
into the theory of normal surface singularities of $(\bc^3,0)$. 
In this way one can use properties of plane curves to get interesting properties of SIS singularities.
They were introduced by I.~Luengo \cite{Ignacio}, and were used 
to answer 
several  questions and conjectures, like the fact that the $\mu$-constant stratum in the semiuniversal deformation space
of an isolated hypersurface singularity is, in general, not smooth. 
Using Zariski pairs as tangent cones of SIS singularities, E.~Artal \cite{ea:jag} found
also counterexamples for a S. S.-T. Yau's conjecture \cite{yau:91} relating the link of the singularity, the characteristic polynomial and the embedded topology.  
A \emph{Zariski pair} is a set of two curves $C_1,C_2\subset\bp^2$ with the same combinatorial type but such that $(\bp^2,C_1)$ is not homeomorphic to $(\bp^2,C_2)$. 
%Of course SIS deserve to be studied also by themselves.

In a recent paper \cite{jag}, A. N\'emethi and the last two authors have found counterexamples 
to several conjectures on normal surface singularities whose link is a rational homology
sphere. For doing this there were used SIS singularities whose tangent cone is a rational cuspidal curve. 
It was shown that the 
\emph{Seiberg-Witten invariant conjecture} (of L.I.~Nicolaescu and A.~N\'emethi \cite{nem:02}), 
the \emph{universal abelian cover 
conjecture} (of W.~Neumann and J.~Wahl \cite{NWnew2})  and the \emph{geometric genus conjecture} (\cite[Question~3.2]{NW}, see also \cite[Problem~9.2]{INV}) 
fail (at least at that 
generality in which they were formulated).

On the other hand, from the positive point of view, SIS singularities have been used by Pi.~Cassou-Nogu\`es and the authors \cite {aclm} 
to confirm the Monodromy Conjecture for the topological zeta function introduced by J.~Denef and F.~Loeser
\cite{dl:92}. We review these results in Section~\ref{sec-mon}.

It is interesting to point out that the relationship between plane  curves and 
normal surface singularities can be used also in the other direction: 
to use results and ideas from normal surface singularities to get new results about curves. 
In this way, 
the results in \cite {aclm} allow to find necessary conditions for the existence of an arrangement of rational plane curves. 
Even more, J.~Fern\'andez de Bobadilla, A. N\'emethi and the last two authors \cite{blmn} 
have found a compatibility property for a rational cuspidal 
projective plane curve to exist based on a heavily study of the failure of the 
Seiberg-Witten invariant conjecture of the corresponding SIS singularities.

Since the class of SIS singularities continue being useful we have decided to write down this survey, dedicated
to our friend Gert-Martin,  where we present 
%a resume of 
known
results and open problems on SIS singularities. 

\section{Superisolated Surface Singularities}\label{sec-sis}

\subsection{Isolated Hypersurface Singularities}

Let $f:(\bc^{n+1},0)\to (\bc,0)$ be an analytic function and the corresponding  germ
$(V,0):=(f^{-1}(0),0)\subset (\bc^{n+1},0)$
of a hypersurface
singularity.
The {\it Milnor fibration}
of the holomorphic function $f$ at $0$ is
the $C^\infty$ locally trivial fibration $f|:B_\varepsilon(0)\cap
f^{-1}(\bd_\eta^*)\to \bd_\eta^*,$ where $B_\varepsilon(0)$ is the open ball
of radius $\varepsilon$ centered at $0$,
$\bd_\eta=\{z\in\bc:|z|<\eta\}$ and $\bd_\eta^*$
is the open punctured disk ($0<\eta\ll\varepsilon$ and $\varepsilon$ small
enough).
Milnor's classical result also shows that the topology
of the germ $(V,0)$ in $(\bc^{n+1},0)$ is determined by the pair $(S_{\varepsilon}^{2n+1},L_V^{2n-1})$,
where $S^{2n+1}=\partial B_\varepsilon(0)$
and $L_V^{2n-1}:=S_{\varepsilon}^{2n+1}\cap V$ is the \emph{link} of the singularity.

Any fiber $F_{f,0}$ of the Milnor fibration is called the  {\it Milnor fiber} of $f$
at $0.$ The {\it monodromy transformation} $h:F_{f,0}\to F_{f,0}$ is the
well-defined (up to isotopy) diffeomorphism of $F_{f,0}$
induced by a small loop around $0\in \bd_\eta.$
The {\it complex algebraic monodromy of $f$ at $0$} is the corresponding linear
transformation $h_*:H_*(F_{f,0},\bc)\to H_*(F_{f,0},\bc)$ on the homology
groups.

If $(V,0)$ defines a germ of isolated hypersurface
singularity then we have that ${\tilde H}_j(F_{f,0},\bc)=0$ but for $j=n.$ In particular
the non-trivial complex algebraic monodromy will be denoted by
$h:H_n(F_{f,0},\bc)\to H_n(F_{f,0},\bc)$ and $\Delta_V(t)$ will denote
its characteristic polynomial. The Monodromy Theorem describes the main properties of the monodromy operator, see for instance the references in  \cite{eb}:
\begin{enumerate}[(a)]
\itemsep0pt
\item $\Delta_V(t)$ is a product of cyclotomic polynomials.
\item Let $N$ be the maximal size of the Jordan blocks of $h$, then $N\leq n+1.$
\item Let $N_1$ be the maximal size of the Jordan blocks of $h$ for the eigenvalue $1$, then $N_1\leq n.$
\end{enumerate}

\subsection{Normal Surface Singularities}

Let $(V,0)=(\{f_1=\ldots=f_m=0\},0)\subset (\bc^N,0)$ be a normal surface singularity with link
$L_V$.
One of the main problems is to determine which analytical properties of $(V,0)$ can be read from the topology of the singularity,
see the very nice survey paper by A.~N\'emethi \cite{INV}. 
Since $V\cap { B}_\varepsilon$ is a cone over the link then $L_V$ characterizes the topological type of $(V,0)$.
%To follow the Artin-Laufer program one needs to restrict to the case of the  link being a rational homology sphere. 

The resolution graph $\Gamma(\pi)$ of a resolution \mbox{$\pi:\tilde V
  \to V$} allows to relate analytical and topological properties of
$V$.  Via plumbing construction, W.~Neumann \cite{neu:81} proved that
the information carried in any resolution graph is the same as the
information carried  by the link $L_V.$ 
Let \mbox{$\pi:\tilde V\to V$} be a \emph{good} resolution of the
singular point $0\in V$. Good means that \mbox{$E=\pi^{-1}(0)$} is a
normal crossing divisor. Let 
$\Gamma(\pi)$ be the dual graph of the resolution (each vertex
decorated with the genus $g(E_i)$ and the self-intersection $E_i^2$ of
$E_i$ in~$\tilde V$). 
Mumford proved that the intersection matrix
\mbox{$I=(E_i\cdot E_j)$} is negative definite and Grauert  
proved the converse, i.e., any such graph corresponds to the link of a
normal surface singularity.

\subsection{Superisolated Surface Singularities}

\begin{definition}
A hypersurface surface singularity $(V,0)\subset (\bc^3,0)$
defined as the zero locus of
$f=f_d+f_{d+1}+\cdots\in\bc\{x,y,z\}$
(where $f_j$ is homogeneous of degree $j$) is \emph{superisolated}, SIS for short,
if the singular points of the complex projective plane curve $C:=\{f_d=0\}\subset \bp^2$ 
% is reduced with
% isolated singularities $\{P_i\}_i$, and these points 
are not situated
on the projective curve $\{f_{d+1}=0\}$, that is $\sing(C)\cap \{f_{d+1}=0\}=\emptyset$ in $\bp^2.$
Note that $C$ must be a reduced curve.
\end{definition}

\noindent
The SIS singularities were introduced by I.~Luengo in \cite{Ignacio}
to study the $\mu$-constant stratum, see Section~\ref{sec-def}.
The main idea is that for a SIS  singularity $(V,0)$, the embedded
topological type (and the equisingular type) of $(V,0)$ does not depend on the
choice of $f_{j}$'s (for $j>d$,
as long as $f_{d+1}$ satisfies the above requirement),
e.g.  one can take $f_j=0$ for any $j>d+1$ and $f_{d+1}=l^{d+1}$
where $l$ is a linear form not vanishing at the singular points \cite{cras}.

\medskip\noindent
{\bf The minimal resolution.} Let $\pi:{\tilde V} \to V$ be the monoidal transformation centered  at the maximal ideal
${\mathfrak m}\subset {\mathcal O}_V$ of the local ring of $V$ at $0.$
Then $(V,0)$ is a SIS singularity if and only if $\tilde V$ is a non-singular surface.
Thus $\pi$ is the {\em  minimal
resolution } of $(V,0)$. To construct the resolution graph $\Gamma(\pi)$
consider $C=D_1+\ldots+D_r$  the decomposition in irreducible components
of the reduced curve $C$ in $\bp^2.$ Let $d_i$ be the degree of the curve $D_i$ in
$\bp^2.$
Then $\pi^{-1}(0)\cong C=D_1+\ldots+D_r$ and the self-intersection
of $D_i$ in $\tilde V$ is $D_i\cdot D_i=-d_i(d-d_i+1)$, \cite[Lemma 3]{Ignacio}.
Since the link $L_V$ can be identified with the boundary of a regular neighbourhood
of $\pi^{-1}(0)$ in $\tilde V$ then the topology of the tangent cone determines the topology of the abstract link $L_V$~\cite{Ignacio}.

\medskip\noindent
{\bf The minimal good resolution} of $(V,0)$ is obtained from $\pi$
by doing the minimal  embedded resolution of each plane curve singularity
$(C,P)\subset  (\bp^2,P)$, $P\in\sing(C)$, which is not an ordinary  double point whose branches belong to different global irreducible components. Let $D_j$ be an irreducible component of $C$ such that $P\in D_j$ and with multiplicity $n\geq 1$ at $P$.
After blowing-up at $P$, the new self-intersection of
the (strict transform of the) curve $D_j$ in the (strict transform of the)
surface $\tilde V$ is $D_j^2-n^2.$
In this way one constructs the minimal good resolution graph
$\Gamma$ of $(V,0)$.

In particular the theory of hypersurface
superisolated surface singularities ``contains''
in a canonical way the theory of complex projective
plane curves.

\begin{example} \label{example} Let $f=f_5+z^6$ be given by the equation
$f_5=z(yz-x^2)^2-2xy^2(yz-x^2)+y^5$.
The curve $C$ is irreducible with unique
singularity at $[0:0:1]$ (of type $\ba_{12}$).
The minimal good resolution graph $\Gamma$ of the superisolated singularity $(V,0)$
is

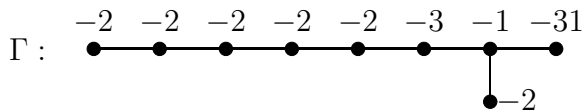
\begin{picture}(300,45)(-30,0)
\put(25,25){\circle*{5}}
\put(50,25){\circle*{5}}
\put(75,25){\circle*{5}}
\put(25,35){\makebox(0,0){$-2$}}
\put(50,35){\makebox(0,0){$-2$}}
\put(75,35){\makebox(0,0){$-2$}}
\put(100,25){\circle*{5}}
\put(125,25){\circle*{5}}
\put(150,25){\circle*{5}}
\put(175,25){\circle*{5}}
\put(200,25){\circle*{5}}
\put(175,5){\circle*{5}}
\put(25,25){\line(1,0){175}}
\put(175,25){\line(0,-1){20}}
\put(0,25){\makebox(0,0){$\Gamma:$}}
\put(100,35){\makebox(0,0){$-2$}}
\put(125,35){\makebox(0,0){$-2$}}
\put(150,35){\makebox(0,0){$-3$}}
\put(175,35){\makebox(0,0){$-1$}}
\put(200,35){\makebox(0,0){$-31$}}
\put(185,5){\makebox(0,0){$-2$}}
\end{picture}
\end{example}

\noindent Here all the curves have genus zero.

\medskip\noindent
{\bf The embedded resolution.} In \cite{ea:mams}, the first author studied the Mixed Hodge
Structure of the cohomology of the Milnor fibre of a SIS singularity. For that he constructed
in an effective way an embedded resolution of a SIS singularity.

The germ  $(V,0)\subset(\bc^{3},0)$ is an isolated surface
singularity. Hence $H_0(F,\bc)$ and $H_2(F,\bc)$ are the
only non-vanishing homology vector spaces on which the
mo\-no\-dro\-my acts (we denote the Minor fiber by $F$).
The only eigenvalue of the action of the monodromy on
$H_0(F,\bc)$ is equal to $1$.
The Jordan form of the complex monodromy  on $H_2(F,\bc)$ is computed for SIS singularity.
Let $\Delta_V(t)$ be the corresponding characteristic polynomial of  the complex monodromy  on $H_2(F,\bc)$. Denote by $\mu(V,0)=\deg(\Delta_V(t))$ the Milnor number of $(V,0)\subset (\bc^3,0)$.

Let $\Delta^P(t)$ be the characteristic polynomial (or Alexander
polynomial)
of the action of the complex
monodromy of the germ $(C,P)$ on $H_1(F_{g^P},\bc),$ (where $g^P$ is a local equation of $C$ at $P$ and $F_{g^P}$ denotes the
corresponding Milnor fiber). Let $\mu^P$ be the Milnor number of $g^P$ at $P$.

Let $H$ be a $\bc$-vector space and $\varphi:H\to H$ an endomorphism of $H$. The $i$-th
Jordan polynomial of $\varphi$, denoted by $\Delta_i(t)$, is the monic polynomial
such that for each $\zeta\in \bc$, the multiplicity of $\zeta$ as a root of $\Delta_i(t)$
is equal to the number of Jordan blocks of size $i+1$ with eigenvalue equal to $\zeta$.

Let $\Delta_1(t)$ and $\Delta_2(t)$ be the first and the second Jordan polynomials of the complex monodromy on $H_2(F,\bc)$ of $V$ and let $\Delta^P_1(t)$ be the first Jordan polynomial of the complex monodromy of the local plane singularity
$(C,P)$. After the Monodromy Theorem these polynomials joint with $\Delta_V(t)$ and $\Delta^P(t)$, $P\in \sing (C)$,
determine the corresponding Jordan form of the complex monodromy.
The Alexander polynomial $\Delta_{C}(t)$ of the projective plane curve $C\subset \bp^2$
was introduced by A.~Libgober \cite{libgober1,libgober2} and F.~Loeser and M.~Vaqui\'e~\cite{lova}.

\begin{theorem}[\cite{ea:mams}]\label{artal} Let $(V,0)$ be a SIS singularity whose tangent cone $C$ 
has $r$ irreducible components. Then  the Jordan form of the complex monodromy  
on $H_2(F,\bc)$ is determined by the following polynomials
\begin{enumerate}[\rm(i)]
\item The characteristic polynomial $\Delta_V(t)$ is equal to
$$
\Delta_V(t)=\frac{(t^d-1)^{\chi(\bp^2\setminus C)}}
{(t-1)}\prod_{P\in\sing(C)}\Delta^P(t^{d+1}).
$$
\item The first Jordan polynomial is equal to
$$
\Delta_1(t)=\frac{1}
{\Delta_{C}(t)(t-1)^{r-1}}\prod_{P\in\sing(C)}\frac{\Delta^P_1(t^{d+1})
\Delta^P_{(d)}(t)}{\Delta^P_{1,(d)}(t)^3},
$$
where
$$\Delta^P_{(d)}(t):=\gcd(\Delta^P(t), (t^d-1)^{\mu^P})\, \text{ and } \, \Delta^P_{1,(d)}(t):=\gcd(\Delta^P_1(t), (t^d-1)^{\mu^P}).$$
\item  The second Jordan polynomial is equal to
$$
\Delta_2(t)=\prod_{P\in\sing(C)}\Delta^P_{1,(d)}(t).
$$
\end{enumerate}
\end{theorem}

\noindent
The first part of the theorem was stated by J. Stevens in \cite{st:89}. A general
formula for the zeta function of the monodromy was proved by D. Siersma \cite{siersma} (see also \cite{zeta},\cite{glm1}).
In particular the Milnor number
$\mu(V,0)$ of a SIS singularity
verifies the identity
\begin{equation*}
\mu(V,0)=(d-1)^3+\sum_{P\in \sing C} \mu^{P}. 
\end{equation*}

\medskip
\noindent {\bf Yomdin Singularities, Series of Singularities and Spectrum.}

\medskip\noindent
The first natural generalization of superislolated singularities
are Yomdin singularities, where $d+1$ is replaced by $d+k$. 

Let
$(V,0)\subset (\bc^{n+1},0)$
be the germ of hypersurface
defined by $f=0$, $f=f_d+f_{d+k}+\ldots\in \bc\{x_0,\ldots,x_n\}$.
The singularity $(V,0)$ is called of \emph{Yomdin type\/}
if \mbox{$\sing(\{f_d=0\})\cap\{f_{d+k}=0\}=\emptyset$} in $\bp^n$. 

For each \mbox{$P\in \sing (\{f_d=0\})$},
let $g^P$ be a local equation of \mbox{$\{f_d=0\}\subset \bp^n$} at $P$.
Formul{\ae} 
for the Milnor number (see \cite{yomdin,le,cras}) and for the zeta function
$\zeta_f(t)$ of the complex monodromy 
can be written as follows \cite{st:89,siersma,glm1,zeta}:
\begin{equation*}
\mu(V,0)=(d-1)^{n+1}+k \sum_{P\in \sing \{f_d=0\}} \mu^{P},
\end{equation*}
and 
\begin{equation*}
\zeta_f(t)=(1-t^d)^{\chi(\bp^{n}\setminus \{f_d=0\})}\left(\prod_{P\in \sing \{f_d=0\}}
(1-t^{d+k}) (\zeta_{g}^{P})^k(t^{d+k})\right)^{-1}.
\end{equation*}

\noindent Here, $(\zeta_{g}^{P})^k(t)$ is the monodromy zeta function
of the $k$-power of the  corresponding monodromy $\zeta_{g}^{P}(t)$ of $g$ at $P.$ 

Let $H$ be a hyperplane such that \mbox{$\sing (\{f_d=0\})\cap H=\emptyset$},
$H$ being the zero locus of a linear form $g$.
Then the family $$F(x_0,\ldots,x_n,t)=f_d+(1-t)(f-f_d)+tg^{d+k}$$
is a $\mu$-constant family (in fact a $\mu^*$-constant family), see \cite{cras}.
It means that to study properties of Yomdin type singularities
which are preserved under $\mu$-constant deformations is equivalent to study \emph{series} 
of singularities of type $f_d+g^{d+k}$. Notice that in such a case the singular locus of $f_d$ is $1$-dimensional.

Let $f$ be a germ of an analytic function at zero whose singular locus
is $1$-dimensional. 
Let $g$ be a generic linear function such that $g(0)=0.$ Y.~Yomdin
\cite{yomdin} compared the  
vanishing cohomologies of their Milnor fibres (and then its Milnor numbers) of
$f$ and $f+g^N$, for $N$ big enough. Later on D.~Siersma  \cite{siersma}
compared the zeta functions of their monodromies. Finally it was J.~Steenbrink
\cite{steenbrink} who conjectured a relationship between the spectrum
$\spec(f,0)$ of $f$ and 
the spectrum $\spec(f+g^N,0)$ of $f+g^N$.
This conjecture was proved by M.~Saito \cite{saito} using his theory of mixed Hodge modules.
Another proof has been given by A.~N\'emethi and J.~Steenbrink, \cite{nem-steen}. Recently
G.~Guibert, F.~Loeser and M.~Merle \cite{glm} have proved Steenbrink's conjecture
without any condition on the singular locus of $f$ and $g$ being any function vanishing at $0.$

The notion of a spectrum $\spec(f,x)$ at $x$ of a function $f$
on a smooth complex algebraic variety was introduced by J.~Steenbrink in \cite{steen1} and by A.~Varchenko in \cite{varchenko}. It is a fractional Laurent polynomial $\sum_{\alpha\in \bq} n_{\alpha} t^\alpha$, $n_\alpha\in \bz$ defined using the semi-simple part of the action of the monodromy on the mixed Hodge structure
on the cohomology of the Milnor fibre of $f$ at $x.$
Here we use the convention given by M.~Saito in \cite{saito} (denoted by $\spec'(f,x)$) which differs from that in \cite{steenbrink} by multiplication by $t$ (see Remark 2.3 in \cite{saito}).

Let \mbox{$f=f_d+f_{d+k}+\ldots\in \bc\{x_0,\ldots,x_n\}$} define a Yomdin type
singularity \mbox{$(V,0)\subset (\bc^{n+1},0)$}
and, for each \mbox{$P\in \sing (\{f_d=0\})$},
let $g^P$ be a local equation of $\{f_d=0\}\subset \bp^n$ at $P.$ Since the
spectrum does not change under $\mu$-constant deformations, see
\cite{varchenko01,varchenko02}, then the spectrum  $\spec(f,0)$ of $(V,0)$ can
be computed via  \cite[Theorem~6.1]{steenbrink} 
and \cite[Theorem~5.7]{saito} in terms of the spectral numbers (also called
exponents) $\{\alpha_i^P\}_P$ of $g^P$ at $P$. 

\begin{theorem}[\cite{cras, saito,steenbrink}]
With the previous notations, the spectrum  $\spec(f,0)$ of a Yomdin singularity $(V,0)\subset (\bc^{n+1},0)$ defined by 
$f=f_d+f_{d+k}+\ldots\in \bc\{x_0,\ldots,x_n\}$ is equal to
\begin{eqnarray*}
\spec(f,0)=&\left(\dfrac{t^{1/d}-t}{1-t^{1/d}}\right)^n
-\left(\dfrac{1-t}{1-t^{1/d}}\right)\sum\limits_{P\in \sing(C)} \sum\limits_{\alpha_i^P \in \spec(g^P,P)}
t^{\beta_i^P}\\
&+\left(\dfrac{1-t}{1-t^{1/(d+k)}}\right)\sum\limits_{P\in \sing(C)} \sum\limits_{\alpha^P_i \in \spec(g^P,P)}
t^{\gamma_i^{P,k}},
\end{eqnarray*}
where $\gamma_i^{P,k}:=\dfrac{k \alpha_i^P+\lfloor d(\alpha_i^P-1)\rfloor +d+1}{d+k}$ and $\beta^P_i:=\dfrac{\lfloor d(\alpha^P_i-1)\rfloor +d+1}{d}$.
\end{theorem}

\noindent
The study of non-isolated singularities defined by an analytic complex function $f$
using perturbation $f+g^k$, $g$ a generic linear form, has been extensively studied
mostly using polar methods, L\^e cycles and other methods, see \cite{massey, massey1,tibar} and references therein.

\section{Deformations}\label{sec-def}

Let \mbox{$p:{\mathcal V} \to {\mathcal T}$} be a deformation of a SIS
singularity \mbox{$(V,0)\subset (\bc^3,0)$} with a section $\sigma$ such that
$(V_t,\sigma(t))$ is an isolated singularity (one may assume that ${\mathcal
  T}$ is one-dimensional 
and smooth). In general $(V_t,\sigma(t))$ is not a SIS singularity but if the
corresponding multiplicities coincides, that is 
\mbox{$\mult(V_t,\sigma(t))=\mult(V,0)$}, then $(V_t,\sigma(t))$ is a SIS
singularity, because one can take local coordinates such that
\mbox{$F(x,y,z,t)=0$} is 
the equation of \mbox{$({\mathcal V},0)\subset (\bc^3\times \bc,0)$},
\mbox{$\sigma(t)=(0,t)$}, 
and $F_d(x,y,z,t)=0$ gives the tangent cone of $(V_t,\sigma(t))$ for all
$t$. Thus $p$ induces a deformation  \mbox{$P:{\mathcal C}\to {\mathcal T}$} of
  the tangent cone \mbox{$C\subset \bp^2$} and since the condition of being a
  SIS singularity is open then \mbox{$\sing(C_t)\cap \{f_{d+1,t}=0)=\emptyset$}
  for $t$ close to $0$. 

Assume now that $p$ is a $\mu$-constant deformation, that is
$\mu(V_t,\sigma(t))=\mu(V,0)$ along the family. Even in this case it is not
known that the multiplicity is constant. In fact the following well-known
problems are still open also for SIS singularities. 

\begin{problem}\label{muct} Let \mbox{$p:{\mathcal V} \to {\mathcal T}$} be a
  deformation of a SIS singularity {$(V,0)\subset (\bc^3,0)$} such that
  \mbox{$\mu(V_t,\sigma(t))=\mu(V,0)$}. Is it true that the multiplicity is
  constant?  
\end{problem}

\begin{problem} Let \mbox{$p:{\mathcal V} \to {\mathcal T}$} be a deformation
  of a SIS singularity {$(V,0)\subset (\bc^3,0)$} such  
that \mbox{$\mu(V_t,\sigma(t))=\mu(V,0)$}. Is it a topologically constant
deformation? 
\end{problem}

\noindent
In \cite[Theorem 1]{Ignacio}, the second author gives an affirmative answer for
Problem~\ref{muct} using B.~Perron results (see \cite{perron}) whose proof
turned out to be incomplete. By putting together \cite[Theorem~2]{Ignacio} and
the correct part of \cite[Theorem~1]{Ignacio}  
%in loc. cit.
one gets:

\begin{theorem} \label{ignacio} Let $p:{\mathcal V} \to {\mathcal T}$ be a
  deformation of a SIS singularity $(V,0)\subset (\bc^3,0)$. Then the following
  conditions are equivalent: 
\begin{enumerate}[\rm(a)]
\itemsep0pt
\item $(V,0)$ is topologically equivalent to $(V_t,\sigma(t))$,
\item $\mu(V_t,\sigma(t))=\mu(V,0)$ and 
$\mult(V_t,\sigma(t))=\mult(V,0)$,
\item the family $\{(V_t,\sigma(t))\}_t$ is $\mu^*$-constant,
\item $\mult(V_t,\sigma(t))=\mult(V,0)$ and the induced deformation
  $P:{\mathcal C}\to {\mathcal T}$ of the tangent cone $C\subset \bp^2$ of
  $(V,0)$ is equisingular. 
\end{enumerate}
\end{theorem}

\noindent
In \cite{Ignacio0}, it was shown how to compute, in an effective way, equations
for the equisingularity stratum $\Sigma_C$ 
of $C$ in the family of all projective plane curves of degree $d$, giving
examples in which $\Sigma_C$ is not smooth. Thus if one considers the SIS
singularity with such tangent cone, then  one gets that  
the $\mu^*$-constant stratum in the versal deformation is not smooth.

The simplest example is \mbox{$f=y(xy^3+z^4)^2+x^9+y^{10}$}. Then $C$ has only
one singular point with an ${\ba}_{35}$ singularity, and  
$\Sigma_C$ is singular. J.~Stevens \cite{st:89} using the $V$-filtration proved
that the $\mu^*$-constant stratum is a component of the $\mu$-constant stratum 
giving the non-smoothness of the $\mu$-constant stratum. 

The nice construction
by V. A.~Vasil'iev and V.V.~Serganova in 
\cite{VS}, using matroids,  gives another examples 
with non-smooth $\mu^*$-constant stratum.
The study of the properties of the equisingularity stratum $\Sigma_C$ of curves
is a classical subject which gained a great impulse with the work of
Gert-Martin Greuel and Ch.~Lossen and E.~Shustin. See \cite{greuel} for a detailed account of
the subject and references. 

Let \mbox{$(V,0)\subset (\bc^3,0)$} be an isolated surface singularity.
B.~Teissier 
asked whether $(V,0)$ can be put in a $\mu^*$-constant family such that there
exists a member of the family which is defined over $\bq$ (resp. $\br$).  
Using SIS singularities one can answer negatively 
to this question. Namely, it is known that there are many curves 
$C\subset \bp^2$ such that no element of the equisingularity stratum
$\Sigma_C$  can be defined over $\bq$ (or $\br$), see \cite{ys:79} and the end
of next section. For such a curve not defined over $\br$ see  
\cite{accm}. If one takes a SIS singularity over such a curve, Theorem
\ref{ignacio} gives us that no member of a $\mu^*$-constant deformation can be
defined over $\bq$ (or $\br$).

\section{Zariski Pairs}\label{sec-zp}

Let us consider $C\subset\bp^2$ a reduced projective curve of degree~$d$ defined by
an equation $f_d(x,y,z)=0$. If $(V,0)\subset(\bc^3,0)$ is a SIS singularity with tangent cone $C$,
then the link $L_V$ of the singularity is completely determined by $C$. Let us recall, that
$L_V$ is a Waldhausen manifold and its plumbing graph is the dual graph of the good minimal resolution. In order to determine $L_V$ we do not need the embedding of $C$ in $\bp^2$, but only its embedding in a regular neighborhood. The needed data can be encoded in a combinatorial way.

\begin{definition}
Let $\irr(C)$ be the set of irreducible components of $C$.
For $P\in\sing(C)$, let $B(P)$ be the set of local irreducible components of $C$. The \emph{combinatorial type} of $C$ is given by:
\begin{itemize}
\itemsep0pt
\item A mapping $\deg:\irr(C)\to\bz$, given by the degrees of the irreducible components of $C$.
\item A mapping $\ltop:\sing(C)\to\Top$, where $\Top$ is the set of topological types of
singular points. The image of a singular point is its topological type.
\item For each $P\in\sing(C)$, a mapping $\beta_P:T(P)\to\irr(C)$
such that if $\gamma$ is a branch of $C$ at $P$, then $\beta_P(\gamma)$ is the global irreducible component containing $\gamma$.
\end{itemize}
\end{definition}

\begin{remark}
There is a natural notion of isomorphism of combinatorial types. It is easily seen that
combinatorial type determines and is determined by any of the following graphs (with vertices decorated with self-intersections):
\begin{itemize}
\itemsep0pt
\item The dual graph of the preimage of $C$ by the minimal resolution of $\sing'(C)$. The set $\sing'(C)$ is obtained from $\sing(C)$ by forgetting ordinary double points whose branches belong to distinct global irreducible components. We need to mark in the graph the $r$ vertices corresponding to $\irr(C)$.

\item The dual graph of the minimal good minimal of $V$.
Since minimal resolution is unique, it is not necessary to mark vertices.
\end{itemize}

\noindent
Note also that the combinatorial type determine the Alexander polynomial $\Delta_V(t)$  of $V$ (see Theorem \ref{artal}).
\end{remark}

\begin{definition}
A \emph{Zariski pair} is a set of two curves $C_1,C_2\subset\bp^2$ with the same combinatorial type but such that $(\bp^2,C_1)$ is not homeomorphic to $(\bp^2,C_2)$.
An \emph{Alexander-Zariski pair} $\{C_1,C_2\}$ is a Zariski pair such that
the Alexander polynomials of $C_1$ and $C_2$ do not coincide.
\end{definition}

\noindent
In \cite{ea:mams},  (see here Theorem \ref{artal}) it is shown that Jordan form of complex monodromy of a SIS singularity is determined
by the combinatorial type and the Alexander polynomial of its tangent cone. 
The first example of Zariski pair was given by Zariski, \cite{zariski,zariski2}; there exist sextic curves with six ordinary cusps. If these
cusps are (resp. not) in a conic then the Alexander polynomial equals $t^2-t+1$ (resp. $1$).
Then, it gives an Alexander-Zariski pairs. Many other examples of Alexander-Zariski pairs have been constructed (Artal,\cite{ea:jag}, Degtyarev \cite{degtyarev}). We state the main result in \cite{ea:mams}.

\begin{theorem}\label{thm-yau}
Let $V_1,V_2$ be two SIS singularities such that their tangent cones form an Alexander-Zariski pair.
Then $V_1$ and $V_2$ have the same abstract topology and characteristic polynomial of the monodromy but not the same embedded topology.
\end{theorem}

\noindent
Recall that the Jordan form of the monodromy is an invariant of the embedded topology of a SIS singularity (see Theorem \ref{artal});
since it depends on the Alexander polynomial $\Delta_C(t)$ of the tangent cone, we deduce this theorem.

\begin{remark}
Every SIS singularity of Theorem~\ref{thm-yau} provides a counter-example to a Conjecture by S.S.T~Yau stated in~\cite{yau:91}:
\emph{abstract topology and characteristic polynomial of the monodromy determine embedded topology}.
\end{remark}

\noindent
There are also examples of Zariski pairs which are not Alexander-Zariski pairs (see \cite{oka,arcar,acct}). Some of them are distinguished by the so-called characteristic varieties
introduced by Libgober~\cite{li:01}. These are subtori of $(\bc^*)^r$, $r:=\#\irr(C)$, which measure the excess of Betti numbers
of finite Abelian coverings of the plane ramified on the curve (as Alexander polynomial does it for cyclic coverings).

\begin{problem}
How can one translate characteristic varieties of a projective curve in terms of invariants of the SIS singularity associated to it?
\end{problem}

\noindent
Though Alexander polynomial and characteristic varieties are topological invariants, they are in fact arithmetic invariants in the following sense. Let us suppose that a curve $C$ is defined by a polynomial with coefficients in a number field $K$; then Alexander polynomial and characteristic varieties can be computed \emph{inside} $K$,
i.e., they do not depend on the embedding $K\hookrightarrow\bc$.

\begin{definition}
An \emph{arithmetic Zariski pair} is a Zariski pair such that its elements are defined with coefficients in a number field and with conjugate equations by the action of a Galois element.
\end{definition}

\noindent
The existence of arithmetic Zariski pairs is a consequence of a work of Serre~\cite{se:64} and Chisini's conjecture~\cite{kulikov}. Explicit examples have been found in~\cite{accm}; moreover, there are a lot of candidates to be arithmetic Zariski pairs,
for example, sextic curves with an $\ba_{19}$ singularity (discovered by Yoshihara~\cite{ys:79}).

\begin{problem}
Let $C_1,C_2$ be an arithmetic Zariski pair and let $V_1,V_2$ SIS singularities such that $C_1,C_2$
are their respective tangent cones. Do they have the same embedded topological type?
\end{problem}

\section{Monodromy Conjecture}\label{sec-mon}
Let \mbox{$f:(\bc^{n+1},0)\to (\bc,0)$} be a germ of a holomorphic function and
let $$(V,0):=(f^{-1}(0),0)\subset (\bc^{n+1},0)$$
be the germ of hypersurface
singularity defined by the zero locus of $f$.

Let $\pi:(Y,\cld)\to(\bc^{n+1},0)$ be
an \emph{embedded resolution} of $(V,0)$, that is, a proper analytic map
on a non-singular complex manifold $Y$ such that:
\begin{enumerate}[(1)]
\itemsep0pt
\item the analytic subspace $\cld:=\pi^{-1}(0)$ of $Y$ is the union of non-singular $n$-dimensional manifolds in $Y$ which are in general position;
\item the map $\pi|_{Y\setminus \cld}$ is an analytic isomorphism: $Y\setminus \cld \to \bc^{n+1}\setminus 0$;
\item in a neighbourhood of any point of $\cld$ there exist a local system of coordinates $y_0,\ldots,y_n$
 such that $f\circ \pi(y_0,\ldots,y_n)=y_0^{N_0}\cdots y_{n}^{N_n}.$
\end{enumerate}
Let
$E_i,i\in I,$ be the irreducible components of the
divisor $\pi^{-1}(f^{-1}(0)).$ For each subset $J\subset I$ we set
$$ E_J:=\bigcap_{j\in J} E_j \,, \ \text{ and } \ \check
E_J:=E_J\setminus\bigcup_{j\notin J} E_{J\cup\{j\}}. $$
For each $j\in I,$ let us denote by $N_j$ the multiplicity of $E_j$ in the divisor
of $f\circ\pi$ and by $\nu_j-1$ the multiplicity of $E_j$
in the divisor of $\pi^*(\omega)$ where $\omega$ is a
non-vanishing holomorphic $(n+1)$-form  in $\bc^{n+1}$.

The invariant we are interested in is
the \emph{local topological zeta function}
$Z_{top,0}(f,s)\in \bq(s),$  which is an analytic
(but not topological, see \cite{aclm:00a}) subtle invariant
associated with any germ of an analytic function $f:(\bc^{n+1},0)\to
(\bc,0)$. This rational function was first introduced by J. Denef and F.
Loeser as a sort of limit
 of the $p$-adic Igusa zeta function, see
\cite{dl:92,dl:98}. The original definition was written in
terms of an embedded resolution of its zero locus germ
$(V,0)\subset (\bc^{n+1},0)$ (although it does not
depend on any particular resolution). 
In
 \cite{dl:98}, J. Denef and F. Loeser
gave an intrinsic definition of $Z_{top,0}(f,s)$
 using arc spaces and the
motivic Igusa zeta function, --\,see also \cite{dl:00a} and the
S\'eminaire Bourbaki talk of E. Looijenga \cite{loo:00}.
\smallbreak

The {\it local topological zeta function} of $f$ is:
\[
Z_{top,0}(f,s):=\sum_{J\subset I}\chi(\check E_J\cap \pi^{-1}(0))\prod_{j\in J}
\dfrac{1}{\nu_j+N_j s}\in\bq(s),
\]
where $\chi$ denotes the Euler-Poincar\'e characteristic.
Each exceptional divisor $E_j$ of an embedded resolution \mbox{$\pi:(Y,\cld)\to
(\bc^{n+1},0)$}  of the germ $(V,0)$ gives a candidate pole $-\nu_j/N_j$ of the
rational 
function $Z_{top,0}(f,s).$ Nevertheless only a few of them
give an actual pole of
 $Z_{top,0}(f,s).$  There are several conjectures related to the topological
zeta functions. We focus our attention  in the {\it Mo\-no\-dro\-my
Conjecture}, see \cite{de:91,dl:92}.

\begin{cjt}[Local Monodromy] If $s_0$
is a pole of the topological zeta function $Z_{top,0}(f,s)$ of the
local singularity defined by $f$, then $\exp(2i\pi s_0)$ is an eigenvalue
of the local monodromy at some complex point of
$f^{-1}(0)$. 
\end{cjt}

\noindent
If $f$ defines an isolated hypersurface singularity,
then $\exp(2i\pi s_0)$ has to be an eigenvalue
of the complex algebraic monodromy of
the germ $(f^{-1}(0),0)$.

There are three general problems to consider when trying to
prove (or disprove) the conjecture using resolution of singularities:
\begin{enumerate}[(i)]
\itemsep0pt
\item\label{notai} Explicit computation of an embedded resolution
of the hypersurface $(V,0) \subset(\bc^{n+1},0).$
\item\label{notaii} Determination of the poles $\{-\nu_j/N_j\}$ of
  $Z_{top,0}(f,s).$ 
\item\label{notaiii} Explicit
computation of the eigenvalues of the complex algebraic monodromy (or
computing the characteristic polynomials of the corresponding action of the
complex algebraic monodromy) in terms of the resolution data.
\end{enumerate}

\noindent
The Monodromy Conjecture, which was first stated for the Igusa zeta function,
has been proved for curve singularities by
F. Loeser \cite{lo:88}. F. Loeser actually proved a stronger
version of the Monodromy Conjecture: \emph{any pole of the
topological zeta function gives a root of the Bernstein
polynomial of the singularity}. The behaviour of the
topological zeta function for germs of
curves is rather well understood once an explicit embedded
resolution \mbox{$\pi:(Y, \cld)\to (\bc^{2},0)$} of curve
singularities is known, e.g.\ the minimal one. Basically,
the poles are the $\{-\nu_j/N_j\}$ coming from rupture 
components in the minimal resolution, see the proof by Veys
\cite{vy:95,vy:97} .  
The case of curves was proved in consecutive works by Strauss, Meuser, Igusa
and Loeser 
for Igusa's local zeta function, but the same proof
works for the topological zeta function.
There are other recent proofs of the conjecture for the case of curves by
Cassou-Nogu\`es and the authors  
\cite{memoirs},  Nicaise \cite{nicaise} and Rodrigues \cite{rodrigues}.

There are other classes of singularities where the embedded resolution  is
known. For example, for any singularity of hypersurface defined by an
analytic function which is non-degenerated with respect to its Newton
polytope,
problems \eqref{notai} and \eqref{notaiii} above are solved. Nevertheless, 
\eqref{notaii} seems to be a hard combinatorial problem. This
problem was partially solved by Loeser in the case where
$f$ has a non-degenerate Newton polytope and satisfies certain
extra technical conditions, --\,\cite{lo:90}.

Even in one of the simplest cases where $f$ has non-isolated singularities,
namely the cone over a curve, problems \eqref{notai} and
\eqref{notaiii} are solved, but  \eqref{notaii} presents serious difficulties. 
B. Rodrigues and W. Veys proved in \cite{rv:pac} the Monodromy Conjecture
for any homogeneous polynomial $f_d\in\bc[x_1,x_2,x_3]$ satisfying
$\chi(\bp^2\setminus \{f_d=0\})\not= 0.$ In \cite{aclm} the authors
complete the proof of this case studying homogeneous polynomials $f_d\in\bc[x_1,x_2,x_3]$ satisfying
$\chi(\bp^2\setminus \{f_d=0\})= 0.$

As we mentioned before,
an embedded resolution is also known for superisolated
surface singularities, --\,see \cite{ea:mams}. This allow Pi.~Cassou-Nogu\`es
and the authors to solve problems \eqref{notaii} and \eqref{notaiii} for SIS
singularities, namely the main result of \cite{aclm} is to prove:

\begin{theorem}[\cite{aclm}]
The local Monodromy Conjecture is true for superisolated surface singularities.
\end{theorem}

\noindent
The
local topological zeta function of a SIS singularity satisfies the following equality,
see \cite[Corollary~1.12]{aclm}:
\begin{multline*}
Z_{top,0}(V,s)
=\frac{\chi(\bp^2\setminus C)}{t-s}+
\frac{\chi(\check C)}{(t-s)(s+1)}+
\\
+\ \sum_{P\in\sing(C)}
\left(\frac{1}{t}+(t+1)\left(\frac{1}{(t-s)(s+1)}-\frac{1}{t}\right)
\!Z_{top,P}(g^P,t)\right),
\end{multline*}
where $g^P$ is a local equation of $C$ at $P$, $\check C:= C\setminus \sing(C)$
and $t:=3+(d+1)s$.

\medskip
\noindent
The following properties can be easily described from the previous equalities:

\begin{proposition}\label{prop-cp} Let $\cp$ be the set of poles of $Z_{top,0}(V,s)$.
\begin{enumerate}[\rm(i)]
\itemsep0pt
\item\label{candpol}
$\cp\subset\{-1,-\frac{3}{d}\}\cup\bigcup
\limits_{P\in\sing(C)}\left\{-\dfrac{3+t_0}{(d+1)}\:\bigg|\:
t_0\text{ pole  of } Z_{top,P}(g^P,t)\right\}$.

\item\label{otros} If $-\frac{3}{d}\neq s_0\in\cp$
then $\exp(2i\pi s_0)$ is an eigenvalue of the complex algebraic monodromy  of $V$.

\item\label{mayor} Let $s_0=-\frac{3}{d}$.
If $s_0$ is a pole of $Z_{top,P}(C,s)$ at some point
$P\in \sing(C)$ and either \mbox{$\chi(\bp^2\setminus C)>0$} or
\mbox{$\chi(\bp^2\setminus C)=0$},
then $\exp(2i\pi s_0)$ is an eigenvalue of the complex algebraic monodromy  of $V$.

\item If $s_0=-\frac{3}{d}$ is a multiple pole of $Z_{top,0}(V,s)$
then $\exp(2i\pi s_0)$ is an eigenvalue of the local algebraic monodromy at some 
singular point of $C$.

\item If $s_0=-\frac{3}{d}$ is not a pole of $Z_{top,P}(C,s)$, the residue of
  $Z_{top,0}(V,s)$ at $-\frac{3}{d}$ equals $d\rho(C)$ 
where
 $$\rho(C):=
\chi(\bp^2\setminus C)+\chi(\check C)\frac{d}{d-3}+
\sum_{P\in\sing(C)}Z_{top,P}\left(C,-\frac{3}{d}\right).
$$
\end{enumerate}
\end{proposition}

\noindent
Following Proposition~\ref{prop-cp}, the Monodromy Conjecture for SIS
singularities is proved in all but two cases: 

\begin{enumerate}[({N}-1)]
\itemsep0pt
\item $\chi(\bp^2\setminus C)= 0$, $s_0=-\frac{3}{d}$ is not a pole
for the local topological zeta function at any singular point in $C$ and
$\rho(C)\neq 0$.

\item $\chi(\bp^2\setminus C)< 0$.
\end{enumerate}

\noindent
The \emph{bad divisors} are the degree $d$ effective divisor $D$
on $\bp^2$ $(d>3)$ such that
$\chi(\bp^2\setminus D)\leq 0$ and $s_0=-\frac{3}{d}$ is not  a pole  of
$Z_{top,P}(g_D^P,s),$ for any singular point $P$ in its support $D_{red}$,
where $g_D^P$ is the local equation of the divisor $D$ at $P.$
The main part of \cite[\S 2]{aclm} is devoted to determine the bad divisors $D$ on
$\bp^2$ such that  $\rho(D) \ne 0$ and finally to prove the Monodromy Conjecture.

Note that the
Euler-Poincar\'e characteristic condition on a bad divisor $D$ implies
that  $D$ has at least two irreducible components, all of them
rational curves, see \cite{djs:95,gp:95,koj:99,alm:ma}.
In fact the main result in \cite{aclm} can be used to study arrangements $C=C_1+\ldots+C_s$
of rational plane curves such that $\chi(\bp^2\setminus C)\leq 0$. In particular
some necessary conditions on the combinatorial type of $C$ (see Section~\ref{sec-zp})
are obtained in order to the curve $C$ exists. 

The authors and S.M.~Gusein-Zade have computed in an unpublished work the topological zeta 
function for Yomdin surface singularities, 
obtaining also a similar formula to the one for SIS singularities.

To avoid problems \eqref{notai} and \eqref{notaii} one can compute the so called motivic Igusa 
zeta function using motivic integration. In particular Pi.~Cassou-Nogu\`es and the authors 
in \cite{memoirs} have verified the conjecture (even the original conjecture by Igusa)
for \emph{quasi-ordinary hypersurface singularities} in arbitrary dimension
measuring arcs and using Newton maps \cite{memoirs}.

\section{SIS Singularities with Rational Homology Sphere Links and
  Rational Cuspidal Curves}

Superisolated surface singularities can be used to construct 
normal surface singularities whose link are rational homology spheres.

Let $(V,0)=(\{f_1=\ldots=f_m=0\},0)\subset (\bc^N,0)$ be a normal surface singularity with link
$L_V$.
One of the main problems is to determine which analytical properties of $(V,0)$ can be read from the topology of the singularity,
see  \cite{INV}. 
Let $\pi:\tilde V \to V$ be a resolution of $V$.

The link $L_V$ is called a \emph{rational homology sphere} ($\bq$HS) if $H_1(L_V,\bq)=\{0\}$,
and $L_V$ is called an \emph{integer homology sphere} ($\bz$HS) if $H_1(L_V,\bz)=\{0\}$. 
In general the first Betti number $b_1(L_V)=b_1(\Gamma(\pi))+2\sum_i g(E_i)$, where 
$b_1(\Gamma(\pi))$ is the number of independent cycles of the graph.
In fact $L_V$ is a $\bq$HS if and only if $\Gamma(\pi)$ is a tree and every $E_i$ is a rational curve. 
If additionally the intersection matrix has determinant $\pm 1$ then $L_V$ is an $\bz$HS.

\begin{example}  If $(V,0)\subset (\bc^3,0)$ is a SIS singularity with an irreducible tangent cone $C\subset \bp^2$ then 
$L_V$ is a rational homology sphere if and only if $C$ is a rational curve and 
each of its singularities $(C,p)$ is locally irreducible, i.e a cusp.
\end{example}

\noindent
In \cite{jag} A.~N\'emethi and the last two authors have used 
SIS singularities whose link is a rational homology sphere to disprove several
conjectures made during last years, see loc. cit. for a series of counterexamples. 
In Example~\ref{example1} we present one of them. 

For instance, it is shown that in the $\bq$HS link case the geometric genus
$p_g$ (analytical property of $(V,0)$) 
does not depend only on its link $L_V$, even if we work only with Gorenstein singularities
(cf. \cite[Question~3.2]{NW}, see also \cite[Problem~9.2]{INV}). 
Moreover for  $\bq$-Gorenstein singularities (with $b_1(L_V)=0$)
analytical properties  like the multiplicity, embedded dimension,
Hilbert-Samuel function are not topological properties. 

It is also shown that the universal abelian cover conjecture by Neumann and
Wahl in \cite{NWnew2} did not held with the generality they stated it. 
The starting point of the conjecture was 
Neumann's result \cite{Neu} that the universal 
abelian cover of a singularity with a good $\bc^*$-action and with
$b_1(L_V)=0$ is a Briesk\"orn complete intersection whose weights 
can be determined from the Seifert invariants of the link.
Their original conjecture was: 

\emph{Assume that $(V,0)$  is $\bq$-Gorenstein 
singularity satisfying $b_1(L_V)=0$. 
Then there exists an equisingular and equivariant 
deformation of the universal  abelian cover of $(V,0)$ 
to an isolated  complete intersection singularity. Moreover, 
the equations of this complete intersection, 
together with the action of $H_1(L_V,\bz)$,
can be recovered from $L_V$ via the ``splice equations''.}

The semigroup condition as stated in \cite{NWnew2} does not hold in general. Thus Neumann and Wahl 
restrict themselves to 
a very interesting class of complete intersection normal complex surface singularities called \emph{splice type singularities},
see \cite{NWnuj1,NWnuj2}. In \cite{NWnuj2} the authors conjectured that rational singularities and $\bq$HS link minimally elliptic singularities belong  to the class of \emph{splice type singularities}. Just recently 
T.~Okuma in \cite{okuma} has given a proof of this result. See the paper by J.~Wahl \cite{wahl} in these proceedings.

Another conjecture that was disproved in \cite{jag}
was the Seiberg-Witten invariant conjecture (SWC). 
A.~N\'emethi  and L.~Nicolaescu \cite{nem:02} offered a candidate as a
topological bound for the geometric genus of a rational homology sphere link of
a normal normal surface singularity:  
\emph{
Let $L_V$ be the link of a normal surface singularity.
\begin{enumerate}[\rm(a)] 
\item\label{swca} If $L_V$ is a rational
homology sphere then 
$$p_g\leq {\bf sw}(L_V)-(Z_K^2+s)/8.$$
\item Additionally, if the singularity is ${\mathbb Q}$-Gorenstein, then
in \eqref{swca} the equality holds.
\end{enumerate}}

\noindent
Here $Z_K$ is the \emph{canonical cycle} associated with $\Gamma(\pi)$, and $s$
the number of vertices in $\Gamma(\pi)$ . Then $Z_K^2+s$ does not depend on the choice of
$\Gamma(\pi)$, it is a topological invariant of $L_V$. Set $H:=H_1(L_V,\bz).$

The {\em Seiberg-Witten invariant} ${\bf sw}(L_V)$ of the link
$L_V$ (associated with the canonical $spin^c$ structure)
is 
$${\bf sw}(L_V):=-\frac{\lambda(L_V)}{|H|}+{\mathcal T}(L_V),$$
where ${\mathcal T}(M)$ is
the sign-refined Reidemeister-Turaev torsion ${\mathcal T}(M)$
(associated with the canonical $spin^c$ structure)  \cite{Tu5} 
and $\lambda(L_V)$ is the normalized by the Casson-Walker invariant, 
using the convention of \cite{Lescop}
(cf. also with \cite{nem:02,[52],[55],INV}). 
Both invariants ${\mathcal T}(L_V)$ and $\lambda(L_V)$ can be determined from the 
graph (for details, see \cite{nem:02} or \cite{INV}).

The SWC-conjecture 
was verified by N\'emethi and Nicolaescu for quotient singularities \cite{nem:02},
for singularities with good 
$\bc^*$-actions \cite{[52]} and hypersurface suspension singularities
$g(u,v)+w^n$ with $g$ irreducible \cite{[55]}.

Let $(V,0)\subset (\bc^3,0)$ be a SIS singularity whose
tangent cone $C\subset \bp^2$ is an  irreducible  rational cuspidal curve
(each singularity of $C$ is
locally irreducible). 

We denote by $\Delta^P$ the characteristic polynomial of
$(C,P)\subset (\bp^2,P)$, set  $\Delta(t) :=\prod\limits_{P\in \sing(C)}\Delta^P(t)$ and
$2\delta:=\deg\Delta(t)$.
By the rationality of $C$ one has   
$$
(d-1)(d-2)=2\delta=\sum_{P\in \sing(C)}\mu^P,
$$
where $\delta $ is the sum of the delta-invariants
of the germs $(C,P)$, $P\in \sing(C)$.

The {\em  minimal
resolution } of $V$ was described in Section~\ref{sec-sis}. 
Since $\Delta^P(1)=1$, this implies that $|H|=\Delta_V(1)=d$.
In fact, one can verify easily that $H=\bz_d$, and a possible generator
of $H$ is an elementary loop in a transversal slice to $C$. 

The other invariants which are involved in the SWC can be computed from the
{\em minimal resolution } of $V$ and using Laufer's formula \cite{laufer}:
\begin{equation}\label{for-sw}
\left\{ \begin{array}{l}
Z_K^2+s=-(d-1)(d^2-3d+1); \ \ \
p_g=d(d-1)(d-2)/6;\ \ \mbox{and} \\ \\
{\displaystyle
{\bf sw}(L_V)=\frac{1}{d}\sum_{\xi^d=1\not=\xi}\ \frac{\Delta(\xi)}{(\xi-1)^2}
 +\frac{1}{2d} \Delta(t)''(1) -\frac{\delta(6\delta-5)}{12d}.
}
\end{array}\right.
%\tag{2}
\end{equation}

\begin{example}\label{example1}
Let us continue with Example~\ref{example}. The link of such SIS singularity 
is a rational homology sphere because the curve $C$ 
is irreducible, rational and cuspidal. The plumbing graph is star-shaped, in particular 
it can be realized by a weighted homogeneous singularity $(V_w,0)$.
 
In this case,
$p_g(V,0)=10$ by the previous formula  and $p_g(V_w,0)=10$ by
Pinkham's formula \cite{Pi1}. In particular, using 
\cite{NW} (3.3), $(V,0)$ is in an equisingular deformation 
of $(V_w,0)$. 
This deformation, found with the help of 
 J.~Stevens, can be described as follows.
The weights of the variables $(a,\ldots,f,\lambda)$ are (62,26,30,28,93,91,-3):
$$V(\lambda)=\left\{ 
\begin{array}{lll}
ab-c^2d=\lambda f\,, & 
bc-d^2=\lambda^2 a\,, & 
ad-c^3=\lambda e\,,  \\
be-df=-\lambda ac^2\,, & 
de-cf=-\lambda a^2\,, & 
af-c^2e=-\lambda b^6\,, \\ 
e^2+a^3+b^6c=0\,, & 
ef+a^2c^2+b^6d=0\,, & 
f^2+ac^4+b^7=0  
\end{array}\right\}\,.$$
Here, \mbox{$(V(0),0)=(V_w,0)\subset (\bc^6,0)$} is Gorenstein, but it is not a
complete intersection. Moreover, 
{\em the two singularities $(V,0)$ and $(V_w,0)$ 
have the same topological types (the same graphs $\Gamma$),
but their embedded dimensions are not the same}: they are 
3 and 6 respectively. It is even more surprising that {\em
their multiplicities are also different}: \mbox{$\mult(V,0)=5$} and
\mbox{$\mult(V_w,0)=6$} (the second computed by {\sc Singular} \cite{GPS01}).  

In \cite{Neu} it was proved that the universal abelian cover $(V_w^{ab},0)$ of $(V_w,0)$ is $\Sigma(13,31,2)$, the Brieskorn 
hypersurface singularity $\{u^{13}+v^{31}+w^2=0\}$.
The corresponding resolution graph $\Gamma^{ab}$ (of both $(V^{ab},0)$
and $(V_w^{ab},0)$) is 

\begin{center}
\begin{picture}(300,45)(0,0)
\put(0,25){\makebox(0,0){$\Gamma^{ab}:$}}
\put(50,25){\circle*{5}}
\put(75,25){\circle*{5}}
\put(100,25){\circle*{5}}
\put(125,25){\circle*{5}}
\put(150,25){\circle*{5}}
\put(175,25){\circle*{5}}
\put(200,25){\circle*{5}}
\put(225,25){\circle*{5}}
\put(250,25){\circle*{5}}
\put(175,5){\circle*{5}}
\put(50,25){\line(1,0){200}}
\put(175,25){\line(0,-1){20}}
\put(50,35){\makebox(0,0){$-7$}}
\put(75,35){\makebox(0,0){$-2$}}
\put(100,35){\makebox(0,0){$-2$}}
\put(125,35){\makebox(0,0){$-2$}}
\put(150,35){\makebox(0,0){$-2$}}
\put(175,35){\makebox(0,0){$-2$}}
\put(200,35){\makebox(0,0){$-2$}}
\put(225,35){\makebox(0,0){$-2$}}
\put(250,35){\makebox(0,0){$-5$}}
\put(185,5){\makebox(0,0){$-2$}}
\end{picture}
\end{center}

\noindent
Even more, \emph{there is no equisingular deformation 
of the universal abelian covers.}
Both $(V^{ab},0)$ and $(V^{ab}_w,0)$ have the same graph $\Gamma^{ab}$
but one can show (see \cite{jag}) that $(V^{ab},0)$
is not in the equisingular deformation of $(V_w^{ab},0)$.

Thus, the only possible  ``splice equation'' 
which defines $(V^{ab},0)$ is $(V_w^{ab},0)$ but 
the universal abelian cover $(V^{ab},0)$ is not in the 
equisingular deformation of $(V_w^{ab},0)$. 
Therefore, the universal abelian cover conjecture is not true. 
Moreover, one has two Gorenstein
singularities (one of them is even a hypersurface Brieskorn singularity)
with the same rational homology sphere link, but with different geometric
genus. This provides counterexample for both SWC 
and geometric genus conjecture. 
\end{example}

\noindent
Looking at the identity \eqref{for-sw}, one considers now the (a priori)
rational function 

\begin{equation}\label{eq-r}
R(t):= \frac{1}{d}\sum_{\xi^d=1}
\frac{\Delta(\xi t)}{(1-\xi t)^2}-\frac{1-t^{d^2}}{(1-t^d)^3}.
%\tag{3}
\end{equation}
J.F.~Fern\'andez de Bobadilla, A.~N\'emethi and the last two authors in \cite{blmn} proved that $R(t)\in \bz[t]$ and it can be written as

\begin{equation}\label{eq-r1}
R(t)=\sum_{l=0}^{d-3} \Big( c_l-\frac{(l+1)(l+2)}{2}\Big)\, t^{d(d-3-l)}\in
\bz[t]. 
%\tag{4}
\end{equation}

\noindent Moreover 
$$R(1)={\bf sw}(L_V)-\frac{Z_K^2+s}{8}-p_g.$$

\noindent In particular, the  ($SWC$) is equivalent to $R(1)\geq 0$.

In fact it is rather curious that in all examples, based on SIS singularities,
studied in  \cite{jag} one gets $R(1)\leq 0.$ 
Motivated by these examples,  in \cite{blmn}
there were worked out many examples discovering that the coefficients of $R(t)$
are always non-positive. 
This gives strong necessary conditions on the singularities of $C$. It is know that in   the problem of classification of rational cuspidal curves one of the key points is to find necessary conditions on the singularities. We state this \emph{compatibility property} on $R(t)$ that we have found as a conjecture.

\begin{cjt}[(CP) \cite{blmn}] 
Let  $(C,p_i)_{i=1}^\nu$ be a
 collection of local plane curve  singularities, all of them locally
irreducible,  such that
$2\delta=(d-1)(d-2)$ for some integer $d$. Then if $(C,p_i)_{i=1}^\nu$
can be realized as the local singularities of a degree $d$
(automatically rational and cuspidal)
projective plane curve of degree $d$ then
\begin{equation}
c_l\leq (l+1)(l+2)/2 \text{ for all \ $l=0,\ldots, d-3$.}
\tag{$*_l$}
\end{equation}
\end{cjt}

\noindent In fact the coefficients $c_l$ can be compute from the polynomial
$Q(t)$ defined in terms of $\Delta(t)$:
$$\Delta(t)=1+(t-1)\delta+(t-1)^2Q(t)=
\sum_{l\nmid d}b_lt^l+\sum_{l=0}^{d-3}c_lt^{(d-3-l)d}.$$

\noindent The main result in \cite{blmn} is to prove

\begin{theorem} If the logarithmic Kodaira dimension
$\bar{\kappa}:=\bar{\kappa}(\bp^2\setminus C)$ is $<2$,
then (CP)
is true. In fact, in all these cases $c_l=\frac{(l+1)(l+2)}{2}$ for any $l=0,\ldots, d-3$. 
\end{theorem}

\begin{corollary}[\cite{blmn}] Let $f=f_d+f_{d+1}+\cdots\ :
(\bc^3,0)\to (\bc,0)$  be a hypersurface superisolated singularity
with $\bar{\kappa}(\bp^2\setminus \{f_d=0\})<2$. Then the Seiberg-Witten
invariant  conjecture  is true for
$(V,0)=(\{f=0\},0)$.
\end{corollary}

\noindent
It is even more interesting to study the compatibility property when $\nu=1$. 
In this case one can prove that 
all the
inequalities ($*_l$) are indeed identities.
These identities are equivalent (via a theorem by A.~Campillo, F.~Delgado and S.M.~Gusein-Zade in \cite{cdg}) to a very remarkable distribution of the elements 
of the semigroup $\Gamma_{(C,P)}$ of the singularity $(C,P)$ in intervals of length $d$.  
It is shown there that $(CP)$ in this case is equivalent to the following conjectural identity:
 
\begin{equation*}
\sum_{k \in \Gamma_{(C,P)}} t^{\lceil k/d\rceil}=
\frac{1-t^d}{(1-t)^2}
=1+2t+\cdots+ (d-1)t^{d-2}+d(t^{d-1}+t^d+t^{d+1}+\cdots).
\end{equation*}

\section{Final Remarks}

\subsection{Weighted-Yomdin Singularities}

The second natural generalization of SIS singularities is obtained if one considers a weighted version of this singularities. 

\begin{definition}
A \emph{weight} is a triple \mbox{$\omega:=(p_x,p_y,p_z)\in\bn^3$} such that
$\gcd(p_x,p_y,p_z)=1$. A polynomial $f$ is $\omega$-weighted-homogeneous of degree~$d$
if $f(t^{p_x} x,t^{p_y} y,t^{p_z} z)=t^d f(x,y,z)$ and defines a curve in the weighted projective plane
\mbox{$\bp^2_\omega:=\bc^3\setminus\{0\}/\sim$},  \mbox{$(x,y,z)\sim
  (t^{p_x},t^{p_y},t^{p_z})$} for all \mbox{$t\in\bc^*$}. 

If \mbox{$P\in\bp^2_\omega$}, we define its order $\nu_P$ as the $\gcd$ of the
weights of the non-zero coordinates of $P$. 
\end{definition}

\begin{definition}
If $C\subset\bp^2_\omega$ is a curve defined by a weighted homogeneous polynomial $f$ and $P\in C$ we define
the \emph{weighted Milnor number} $\mu^\omega(C,P)$ as $\frac{\mu}{\nu_P}$ where $\mu$ is defined as follows; let us suppose that $P$ is the equivalence class of $(x_0,y_0,1)$ and consider the Milnor number of $f(x,y,1)=0$
at $(x_0,y_0)$. A \emph{singular point} of $C$ is a point such that $\mu^\omega(C,P)>0$.
\end{definition}

\noindent
Let us consider a germ \mbox{$(W,0)\subset (\bc^3,0)$} defined by a power
series $g$; let 
$g=g_d+g_{d+k}+\dots$ be the weighted-homogeneous decomposition of $f$ with
respect to $\omega$ 
and let $C_m^\omega\subset\bp^2_\omega$ be the weighted-projective locus of zeroes of $g_m$. 

\begin{definition}
We say that $(W,0)\subset (\bc^3,0)$ is a \emph{weighted-Yomdin singularity with respect to $\omega$} if $\sing(C^\omega)\cap C_{d+k}^\omega=\emptyset$.
\end{definition}

\noindent
In a forthcoming joint work with F. Fern\'andez de Bobadilla, we will give a proof of a formula which was suggested to us by C.~Hertling.
\begin{proposition}
The Milnor number $\mu$ of a \emph{weighted-Yomdin singularity $(W,0)\subset (\bc^3,0)$ with respect to $\omega$} satisfies the following equality:
$$
\mu(W,0)=\left(\frac{d}{p_x}-1\right)\left(\frac{d}{p_y}-1\right)\left(\frac{d}{p_z}-1\right)+
k\sum_{P\in\sing(C^\omega)}\mu(C^\omega,P).
$$
\end{proposition}

\subsection{$*$-Polynomials}

The theory of (local) SIS or Yomdin singularities has an analogous global counterpart defined by polynomials
of type $f=f_d+f_{d-k}+\dots$ and the same geometric condition.
For instance the formula for the global Milnor number is done by the authors in 
\cite{alm} and for the zeta-function of the monodromy at infinity by S.M.~Gusein-Zade and the last two authors in 
\cite{glm2}. 
A finer study has been done in a series of works A.~N\'emethi and R.~Garc\'{\i}a L\'opez \cite{gn1,gn2,gn3} for $*$-polynomials $f=f_d+f_{d-1}+\dots$. The behaviour of these polynomials at infinity imitates in some way the local behaviour of SIS singularities. They computed formul{\ae} for the global Milnor number, monodromy at infinity, Mixed Hodge structure at infinity...

\subsection{Intersection form of a SIS singularity}
 
In the topological study of singularities, we are interested in invariants living in the complex setting (like the Jordan form of the monodromy) but also in invariants living in the integers, like monodromy over $\bz$, Seifert form
or the intersection form in a distinguished basis of vanishing cycles.

It is well-known how to compute these invariants for local germs of curves. In his thesis, M.~Escario computes these invariants for polynomials in two variables which are generic at infinity (in fact, for the more general concept of \emph{tame} polynomials), using a generic polar mapping $\Phi$ and the braid monodromy of the discriminant of $\Phi$.

Combining these techniques with Gabri\'elov's method (see \cite{gabrielov}), M.~Escario gives a method to compute the intersection form of the Milnor fiber in a distinguished basis of vanishing cycles for SIS singularities. In fact, this method works also for Yomdin singularities. 

\subsection{Durfee's conjecture for SIS singularities}

A.~Durfee \cite{durfee} conjectured that the signature of the Milnor fibre of an hypersurface surface singularity  is negative.  In fact, Durfee's conjecture is the stronger inequality
\begin{equation*}
6p_g\leq \mu. \tag{*}
\end{equation*}
Y.~Xu and S.S.T.~Yau proved $(*)$ for weighted-homogeneous surface singularties, \cite{xy}. A.~Nem\'ethi
\cite{nem:98} verified $(*)$ in  the case $f(x,y)+z^n$ with $f(x,y)\in \bc\{x,y\}$ irreducible, see also 
T.~Ashikaga \cite{ashikaga}.

Using SIS singularities,  A. Melle Hern\'andez \cite{israel}
proved $(*)$ for \emph{absolutely isolated} surface singularities.
A surface hypersurface $(V,0)\subset (\bc^3,0)$ is absolutely isolated if there exists
a resolution $\pi:\tilde V \to V$ such that $\pi$ is a composition of blowing-ups at points.

\section*{Acknowledgments}

The authors would like to thank the persons who have collaborated with them in their
works on SIS singularities  whose ideas and comments
are cited in this survey, in special to  
Pi.~Cassou-Nogu\`es, J.~Carmona, J.I.~Cogolludo, M.~Escario,
J.~Fern\'andez de Bobadilla, S.M.~Gusein-Zade, C.~Hertling, 
M.~Marco, A.~N\'emethi, J.~Stevens, H.~Tokunaga and W.~Veys.

\noindent
{\sc Departamento de Matem\'aticas, Universidad de Zaragoza,
Campus Pza.\ San Francisco s/n,
E-50009 Zaragoza, Spain}\\
{\it E-mail address:} \texttt{artal@unizar.es}

\medskip\noindent
{\sc Departamento de \'Algebra,
Universidad Complutense,
Plaza de Ciencias 3,
E-28040 Madrid, Spain}\\
{\it E-mail address:} \texttt{iluengo@mat.ucm.es}

\medskip\noindent
{\sc Departamento de \'Algebra,
Universidad Complutense,
Plaza de Ciencias 3,
E-28040 Madrid, Spain}\\
{\it E-mail address:} \texttt{amelle@mat.ucm.es}

%\newpage
%\tableofcontents

\end{document}